# Algorithm for 2D Mesh Decomposition in Distributed Design Optimization

Shuvodeep De[a,b,*], and Rakesh K Kapania[a,*]

[a]Kevin T. Crofton Department of Aerospace and Ocean Engineering, Virginia Polytechnic Institute and State University, Blacksburg, VA, 24061-0203, USA; [b]Department of Chemical and Biological Engineering University of Alabama, Tuscaloosa, AL, 35487, USA

**ABSTRACT**

Optimization of thin-walled structures like an aircraft wing often involves dividing it into numerous localized panels, each characterized by its own set of design variables. In this work, a generalized algorithm to extract localized panels from the two-dimensional (2D) mesh is discussed. The process employs set operations on elemental connectivity information and is independent of nodal coordinates. Thus, it is capable of extracting a panel of any shape given the boundary and can be used during the optimization of a wide range of structures. A method to create stiffeners on the resulting local panels is also developed. The panel extraction process is demonstrated by integrating it into a distributed MDO framework for optimization of an aircraft wing having curvilinear spars and ribs (SpaRibs). A range of examples is included wherein the process is used to create panels on the wing-skin, bounded by adjacent SpaRibs.

**INTRODUCTION**

One of the most important concerns while designing vehicles is to reduce structural weight, which can directly lead to a reduction in fuel consumption. During the era when computers were expensive and not so powerful, structural design was mostly done by hand calculation on simplified mathematical models. However, since the seventies, due to the rapid increase in computational power numerical solution techniques like Finite Element Analysis (FEA) have gained immense popularity. Not only the details and the complexity of the system can be included in the analysis, but multiple disciplines can now be considered while setting up a structural optimization problem. This avenue of research where several disciplines are incorporated into the optimization problem is known as Multidisciplinary Design Optimization (MDO). The major advantage of solving such a problem arises when relevant disciplines are not independent of each other; in other words, the disciplines interact with each other.

Application of MDO for structural design can be traced back to the work of Schmit [1,2,3]. In these work, finite element methods and algorithms for numerical optimization were used. In subsequent years, Starnes, J. and Haftka[4], Haftka et al. [5] and Fulton et al. [6] designed aircraft wings considering constraints on strength, stability and flutter velocity. With the help of commercial softwares like ANSYS, MDO can be applied to complex 3D spacial structures as well (For example as done by Aru et al.[7]). With the advance in computational power, MDO rapidly gained popularity power in aerospace (Li et al. [8]) and automobile engineering and soon, problems were solved involving a complete three dimensional structures (Devarajan et al. [9]; Singh et al. [10]; Fulton et al. [11]; Kroo et al. [12], Manning [13]; Shihan [14]; Antoine and Kroo[15]; Henderson et al.[16]; Alonso and Colonno[17]). The processes followed in MDO have either monolithic or distributed architectures (Martins and Lambe[18]). In any MDO, the first step is almost always to describe the system using a set of design variables which are often selected by Design of Experiment (DOE) techniques (An example of DOE can be found in the article by Kumar et al. [19]). The goal is to find the best values for the design variables that minimizes (or maximizes) the objective function while satisfying constraints in several disciplines. For example, in problems involving structural design, the size, shape or topology of the structure are described by a set of design variables. By applying appropriate numerical optimization algorithms, the problem is solved for a set of design variables that gives minimum weight or maximum compliance while satisfying constraint like maximum von Mises stress, minimum buckling factor, maximum displacement etc. In a monolithic architecture, all the design variables and constraints are considered in a single optimization process. This process, commonly known as All-at-Once optimization (Sobieszczanski-Sobieski and Haftka[20]), although simple to implement, is computationally expensive when the number of design variables is large. Such problems involving a large number of design variables are often solved using the other process of MDO i.e. distributed architecture. The distributed MDO architectures involves the decomposition of complex systems into multiple smaller components which are then described by a lower number of design variables and optimized independently. The process is usually implemented using parallel computation which can reduce wall clock time by several times. This process of decomposition of a system into a simple sub-system is often known as global/local design optimization.

The aircraft wing is a complex structure consisting of the outer aerofoil shell known as the wing-skin, and internal stiffening elements: the ribs and spars. The global/local optimization of wing design has been used by several research groups including Cimpa et al. [21], Yang et al. [22] and Barkanov et al. [23]. Even though the availability of computational power makes the exploration of large design space feasible, there has always been a concern in the industry about manufacturing limitations. It is often very expensive to produce unconventional designs using conventional manufacturing processed. However, with the invention of 3D printing techniques, the manufacturing industry is likely to be revolutionized over the next few decades. A new additive manufacturing technique known as Electron Beam Free Form Fabrication or EBF3, in short, has recently been developed by Taminger and Hafley[24] at NASA Langley Research Center to fabricate metallic structures of complex shapes, which now can be printed with significant precision. This technology inspired Kapania et al. at Virginia Tech (Locatelli et al. [25]; Locatelli et al. [26]; Jrad et al.[27]; Devarajan et al.[28]; Miglani et al.[29]) to propose the use of curvilinear stiffening elements to reduce structural weight and achieve desirable aeroelastic properties for aircraft. The EBF3GLWingOpt is one of the several optimization frameworks that is being developed at Virginia Tech to optimize aircraft structures using curvilinear spars and ribs (*SpaRibs*). It performs global/local optimization of high aspect ratio cantilever transport aircraft wing for multiple constraints including stress, buckling and crippling. The wing geometry and mesh are generated using commercial software, MSC.PATRAN and MSC.NASTRAN is used for static and buckling analysis. The framework is written in *python* environment and it enables the use of parallel processing.

In the original version of EBF3GLWingOpt, written by Liu et al., the Linked-shape method proposed by Locatelli et al. [27] is used to create *SpaRibs* in each of the wing-boxes using a limited number of design variables (which specified the shape of the *SpaRibs*). The upper and the lower skin of the wing are divided into local panels using the intersection of the *SpaRibs* and the stiffeners are attached to each of the panels. The thickness of each of the local panels, the stiffener height and thickness are considered as design variables. The optimization framework could only explore limited design space where *SpaRibs* could start at the leading-edge spar and end at the trailing-edge spar. In order to overcome this limitation, De et al. [30] proposed the Extended-space method to create *SpaRibs*. By this method not only the constraint on the starting and ending point of the *SpaRibs* was removed but also *SpaRibs* crossing the junction of the inner and outer wing-box can be created. In addition, the algorithm to extract local panel from the finite element model as used the original version of the EBF3GLWingOpt framework was dependent on nodal coordinate coordinates of the nodes and could extract only panels with four edges from the finite element model. Moreover, the idea of dividing a surface into local panels and assigning a thickness design variable to each of the panels is very generalized and can be applied not only to aircraft wings but different other structures including automobile, ships, buildings etc. Thus, the need was felt to develop a generalized algorithm to break a finite element model into local panels. In the following work, such an algorithm to divide a surface mesh (consisting of triangular elements) has been discussed. The process can be implemented on the CRM wing for any *SpaRibs* configuration as well as other shell-structures like the fuselage or automobile frame with minor modifications. This algorithm is not dependent on nodal coordinates and is purely based on set operations performed on the element connectivity matrix.

The article is organized intothreeparts. First, a detail description of the developed algorithm to extract the local panel is given. Second,the capability of the methodis described using different examples. In the first example, the method is implemented on a simple rectangular wing-box to create local panels from the top and bottom skin. In the second example, the method is integrated with the EBF3GLWingOpt framework and demonstrates its effectiveness for a range of *SpaRibs* profiles with examples. Finally, a method of laying down stiffeners on the panels after they are generated using the algorithm and performing optimizationis demonstratedbefore the concluding article.

**METHODOLOGY**

The goal is to develop an algorithm to split a 2D finite element surface mesh into local panels given the boundary nodes of the panel. To do so, the elements interior to the local panel along its outer edge, i.e., the element which shares the boundary nodes, need to be determined. Once these boundary nodes are identified, all nodes and elements interior to the local panel can be found using the connectivity matrix information, following the process mentioned in the flowchart given in Figure 1. Information about the nodal coordinates is not required for this process. Figure 2 shows a general 2D surface mesh from which local panels are needed to be generated. Figure 2a shows the nodes (marked with black dots) are the boundary of the panel. Figure 2b shows, in orange, the boundary elements and interior to the panel. Using the process which is described in Figure 1, all elements of the panel (shown in orange in

Figure2c) are determined. The list of elements and list of nodes of the local panel are stored in the lists: ELEMENT_LIST and NODE_LIST, respectively.

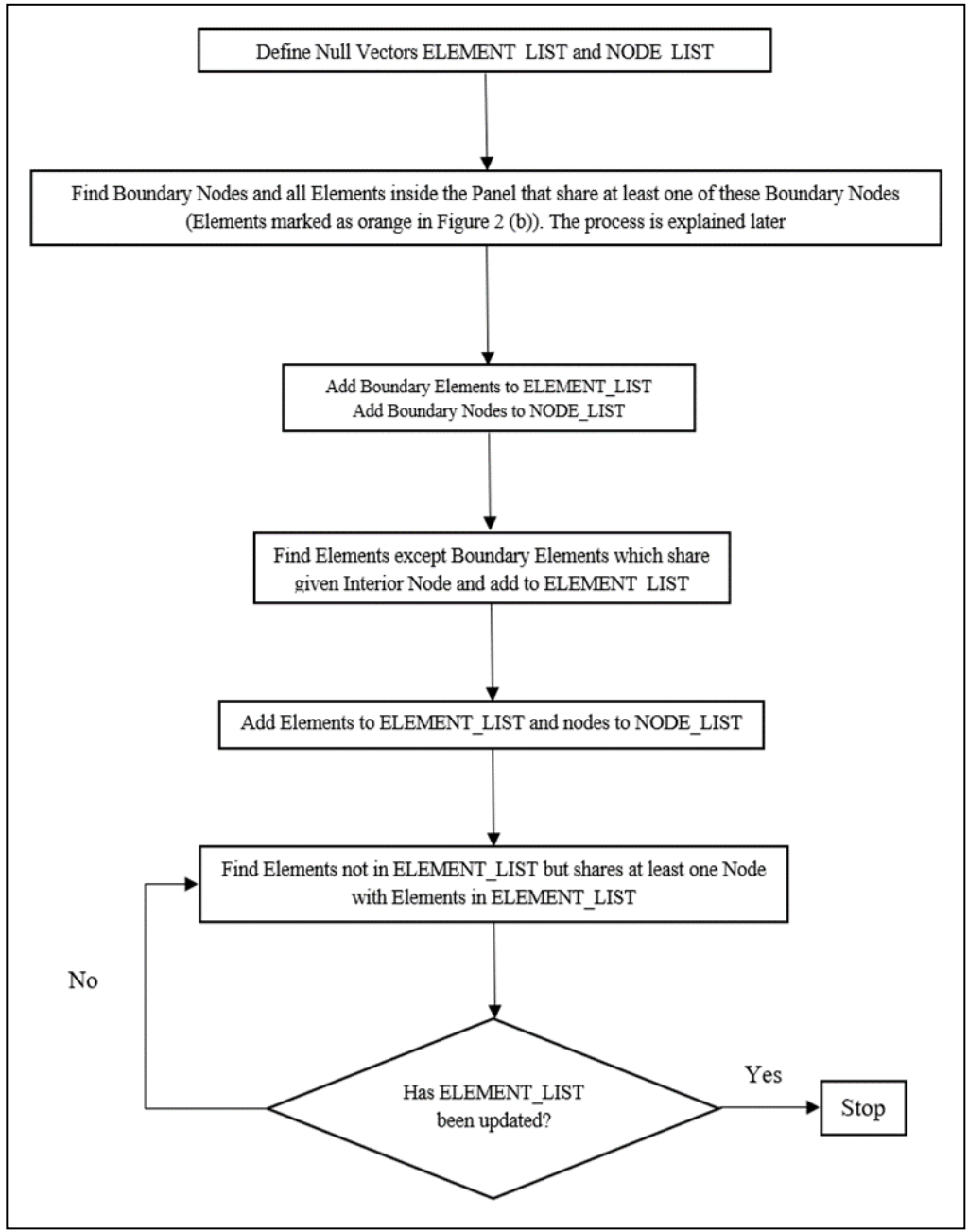

**Figure 1. Algorithm to find Local Panel (De et al. 2018)**

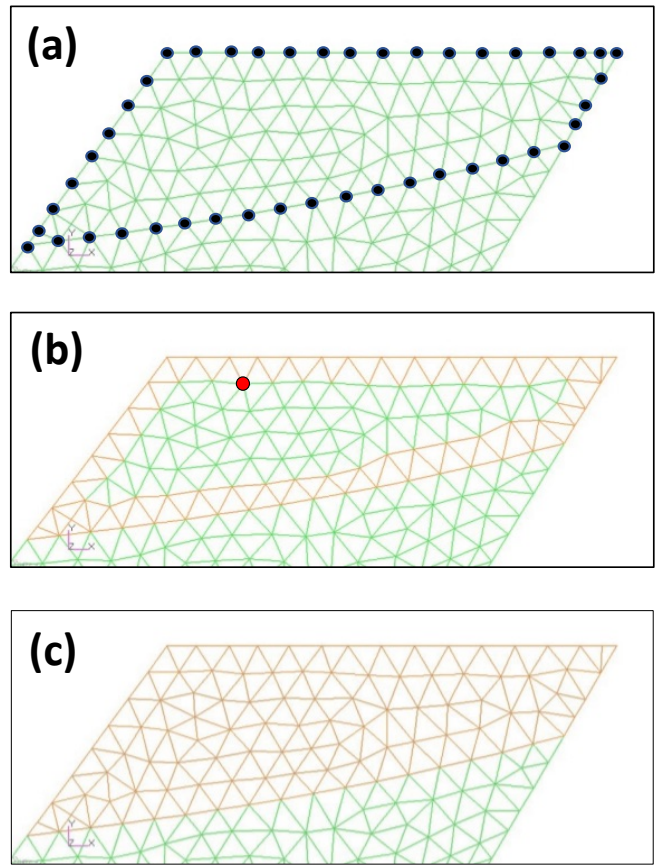

**Figure 2. Determination of elements constituting local panel: (a) Boundary nodes of panel marked with black dots; (b) Boundary elements and interior starting node (marked with a red dot); (c) Elements comprising the local panel.**

### *The algorithm todetermine the boundary elements*

The Mesh-continuity algorithm is easy to apply once the elements along the outer boundary of the panel are determined. The general process to find these boundary elements is complicated and will be discussed in detail in this section with a simple example where the objective is to split a rectangular plate (with holes) that has been meshed with triangular elements, into two panels. *Since the process does not depend on the nodal coordinates, it would work on any mesh created by alinear transformation of this flat-plate.* The process consists of two steps. The first step is to determine all elements along the outer periphery of the panels. The second step is to determine the elements with edges along the 'dividing curve,' which lies inside the finite element mesh. In Figure3, the outer periphery of the panel and dividing curve are shown.

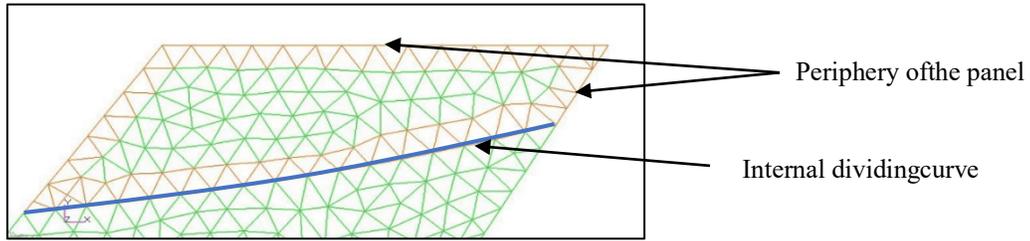

**Figure 3.** The periphery of the panel and the internal dividing curve.

### *The algorithm to determine elements outer periphery*

The plate shown in Figure4 is needed to be split along Curve 2 (Curve 3 is the next curve in the family; Curve 1 is the previous one). The first elements along the periphery of the panels, will be found. In the algorithms described in this section, CHECKPOINT_NODE and CHECKPOINT_ELEMENTrefer to the node and element respectivelywith respect to which the positions of the next node determined. ELMENT_LIST is a list containinglists of the form: [Element reference number, Connectivity Nodes] while NODE_LIST isa list ofthe reference number of the nodes. All of them are initialized as null lists. The algorithm is described in Figure5 using a flowchart.

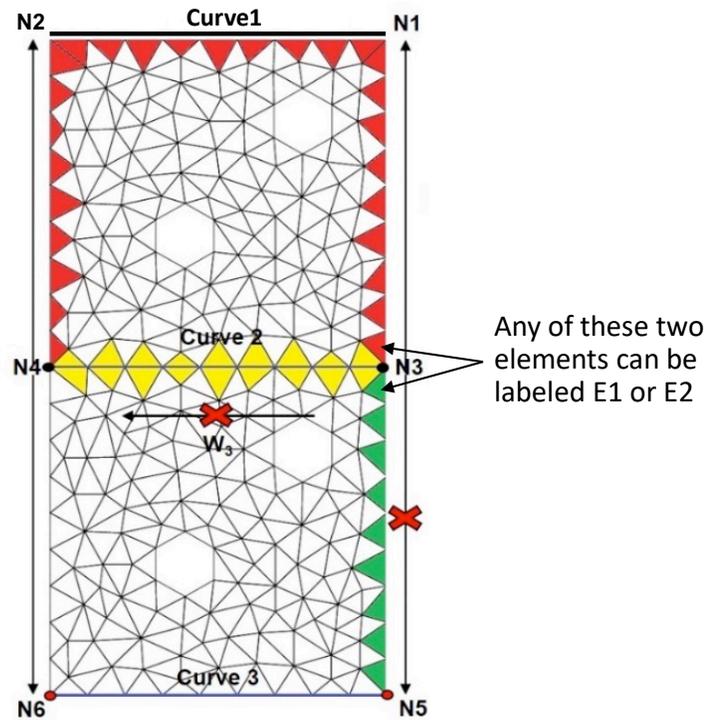

**Figure 4. Application to splitting a plate.**

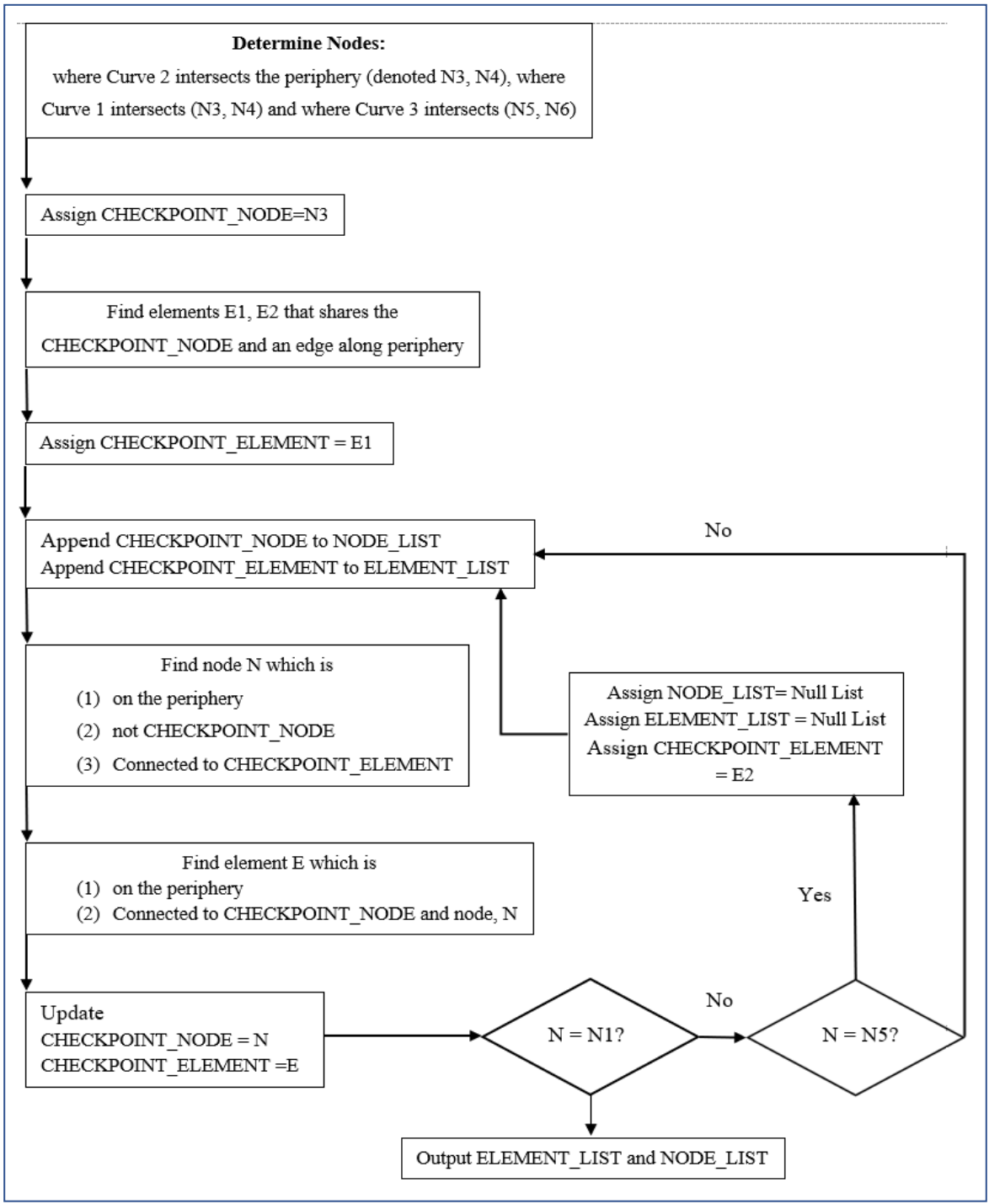

Figure 5. Algorithm to find elements along the outer boundary.

The nodes along Curve1 are already known. Thus, finding the elements along this curve is straightforward. All it needs is to find the elements that have at least two nodes along Curve1. The elements in ELEMENT_LIST and those along Curve1 are marked in red in Figure4. If the CHECKPOINT_NODE is assigned N5, it means we are moving along the elements marked in green, which is not the correct path. It also means that we have started with the element among E1 and E2 which lies outside the panel. If we move in the wrong path, it requires the re-initialization of ELEMENT_LIST and NODE_LIST as a null list and starts the process again starting from the correct element among E1 and E2.The process described above cannot be used to find elements interior to the local panel with edges along Curve2 as it will include elements outside the local panel as well (elements shown in yellow in Figure4), because all these elements have edges along Curve2. To find the elements located inside the panels, another algorithmknown as Mid-element Algorithm (MEA) has been developed, which is described in the flowchart in Figure6.

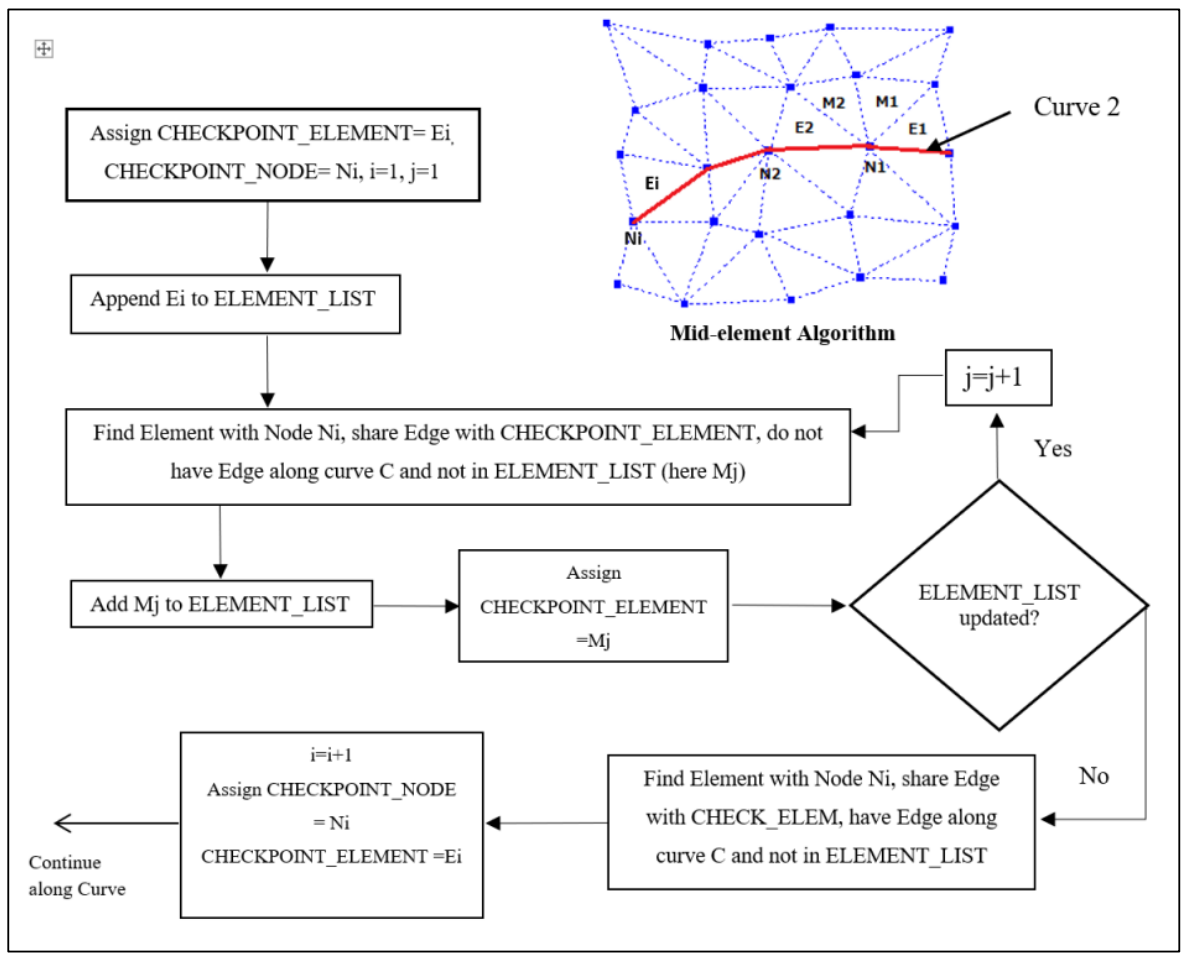

**Figure 6. Flowsheet diagram showing the middle element algorithm (MEA).**

While implementing theMid-element algorithm (MEA), we collect all the elements and nodes in the same lists, i.e.,ELMENT_LIST and NODE_LIST, respectively.Similar to the process used to determine elements along the periphery, we use two variables CHECKPOINT_NODE and CHECKPOINT_ELEMENT for node and element with respect to which the next node and element are found. The algorithm is applied along Curve 2 considering the nodes in order. The reference number of the nodes and the element interior to the panel along the curve, Curve 2, in order from right to left, are $N_i$ and $E_i$, respectively. The stopping criteria are when CHECKPOINT_NODE = N4.

However, before implementing the MEA, it is important to correctly choose the first element i.e. E1. It must be an element interior to the panel. This can be determined with the help of the current CHECKPOINT_ELEMENT, which contains only the elements marked with red in Figure 7. We initiate another vector, CHECKPOINT_ELEMENT_1 and a null list, ELEMENT_LIST_1. The process is given in Figure 8. Here the 'Mid-element algorithm' (MEA) is applied to the flat plate example:

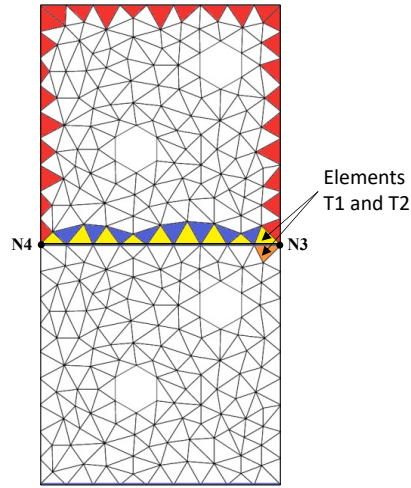

**Figure 7. Using the Middle-element algorithm to find elements colored yellow and blue in Local panel.**

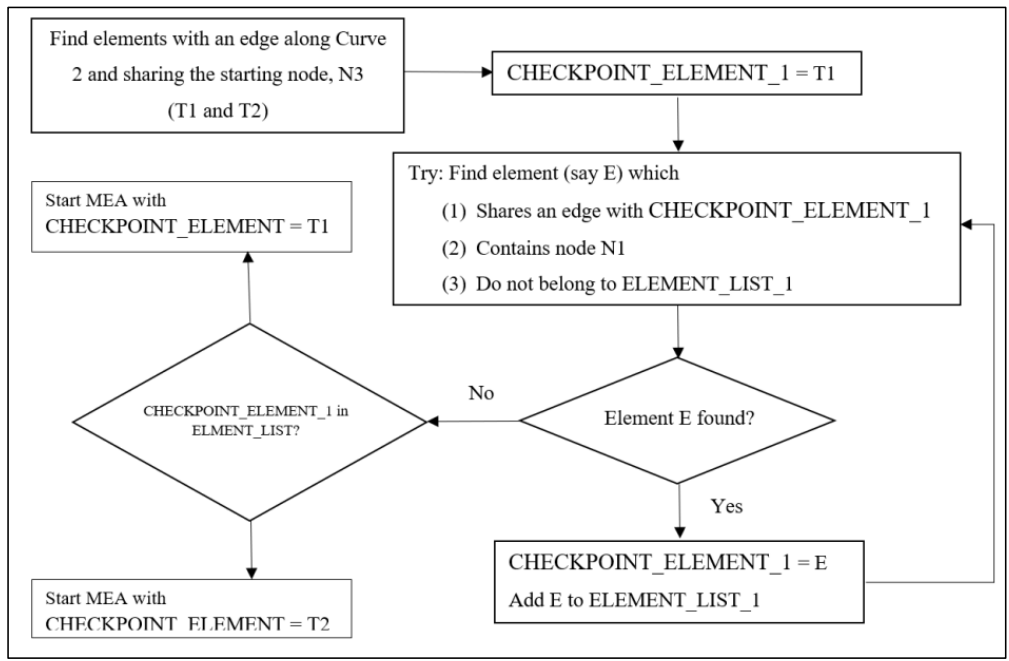

**Figure 8. Algorithm to find the correct first element for MEA.**

The elements included in the list ELEMENT_LIST after implementation of MEAinclude elements above Curve 2 (marked by yellow and blue in Fig.8 and elements form the periphery of the panel above Curve 2 (marked in red in Figure8). Once these elements along the outer periphery are determined, the Mesh-continuity algorithm described at the start of this section can be used to find all the elements belonging to the panel.

**RESULTS AND DISCUSSION**

*Demonstration on rectangular wing-box*

The method is demonstrated first on a simple example. The rear wing-box of Boeing HSCT N+2 Wing containing 7 *SpaRibs* in the span-wise direction and 8 *SpaRibs* in the chord-wise direction is constructed in MSC.PATRAN, as shown in Figure9a. Figure 9b shows a local panel extracted from the top skin from the wing-box using the method.

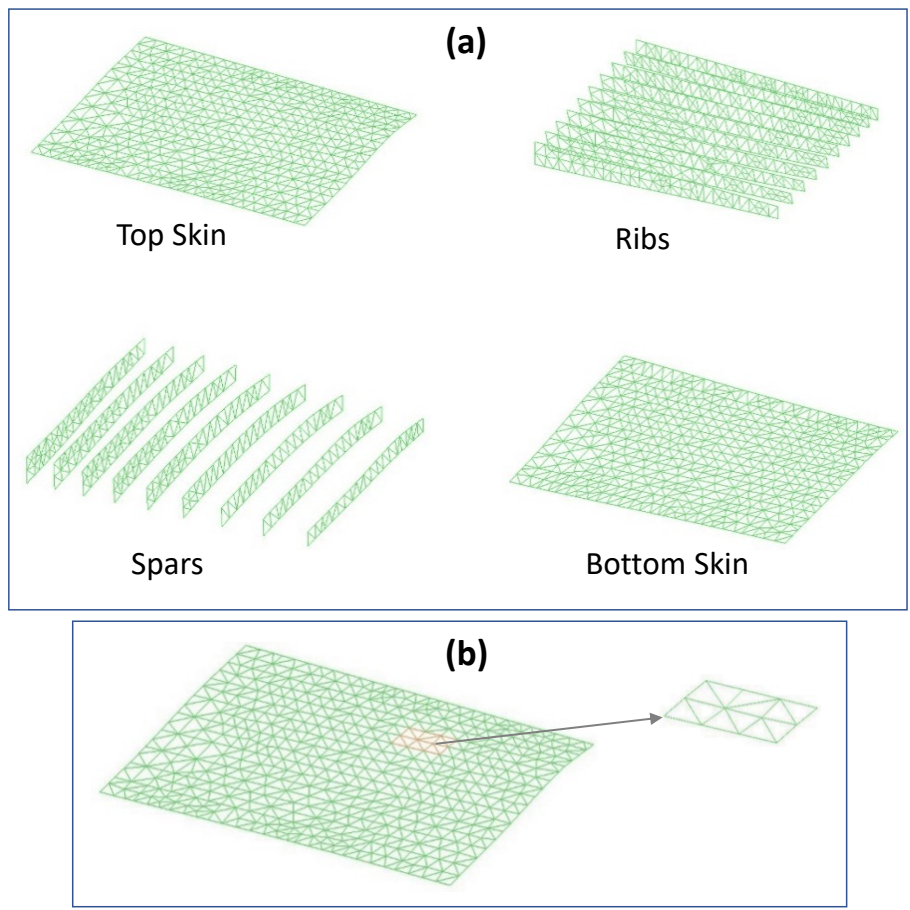

Figure 9. (a)Mesh of rear wing-box of Boeing HSCT N+2 wing using triangular elements; (b) Local panel from upper skin of Boeing N+2 HSCT wing-box.

*Integration of algorithm with EBF3GLWingOpt*

The algorithm can be integrated with any distributed MDO framework if the purpose is to optimize the thickness of 2D surfaces. EBF3GLWingOpt is one such MDO framework being developed at Virginia Tech by Kapania et al. to optimize commercial aircraft wing with curvilinear spars and ribs (*SpaRibs*). The shape of the *SpaRibs*, as well as thickness distribution of the wing-skin and the *SpaRibs*, are considered for optimization.

The algorithm is integrated with the EBF3GLWingOpt framework to create a local panel on the upper and lower skin of the NASA CRM Wing (high-aspect ratio transport aircraft wing) for a range of *SpaRibs* configurations possible using the Extended-space Method. The wing-skin is first divided by the set of *SpaRibs,* which replace the ribs. Each of the sections is then divided by the family of *SpaRibs* that replace the spars. The boundary nodes are the nodes common to the wing-skin and the *SpaRibs*, thus bounding the local panel. Since the algorithm is completely based on set operations performed on the connectivity matrix and independent of nodal coordinates, the algorithm can be used in a variety of problems requiring the generation of local panels. The process is also independent of the element reference numbers, and it works no matter in whatever order the elements are distributed. Depending on the shape of the bounding SpaRibs, the panel can have any number of edges greater than three. Another advantage of this method is that the boundary nodes of each of the local panels are already determined, and the information is stored to be used in imposing boundary conditions during the optimization process. The proposed algorithm has been implemented on wings with different SpaRibs profiles and is found to be useful for SpaRibs of different orientations and curvatures. Figure 10a shows the internal spars and ribs layout of a commercial aircraft wings. The wing contains one inner spar and 37 straight ribs. Figures10b and c show the top skin of the wing which is divided into local panels using the algorithm and some of the local panels created (marked in orange in b), respectively.Figure 10d and e show the internal layout of a wing with curvilinear spars and ribs (SpaRibs) of two different profiles, respectively. The local panels created from top skin using SpaRibs shown in Figure 10d and e are represented in Figure 10f,g and Figure10h,i, respectively.

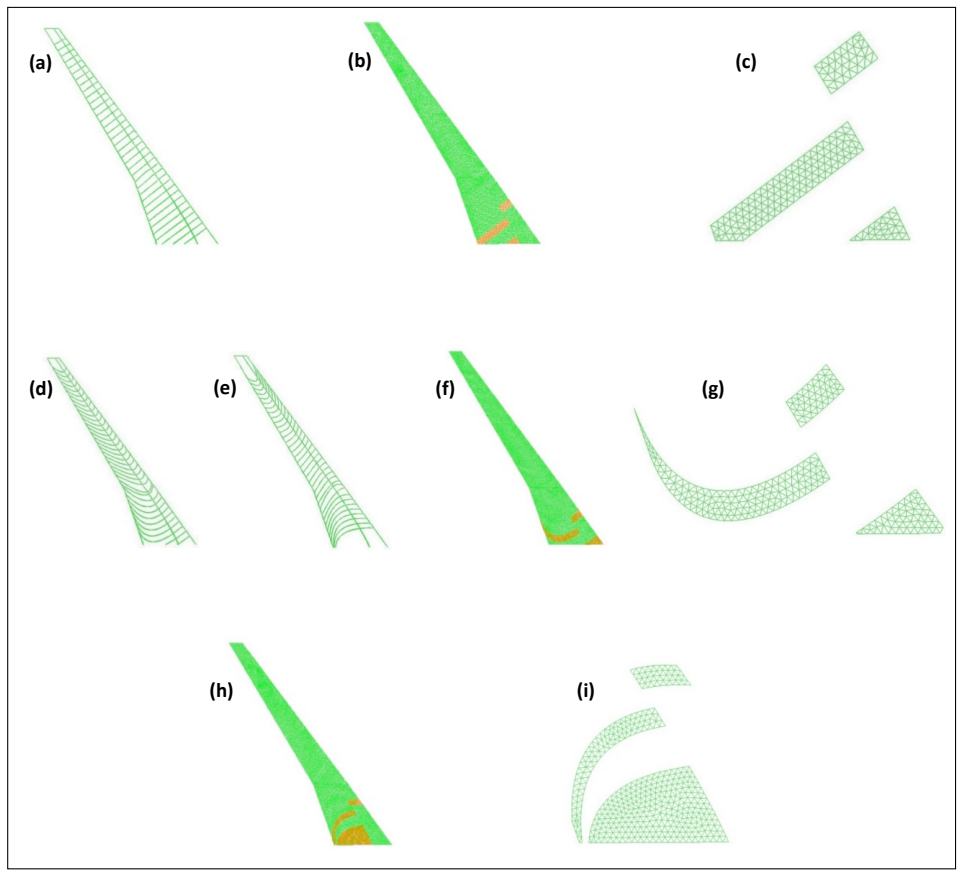

Figure 10. (a) Internal layout of a commercial aircraft wing (with one internal Spar and 37 straight Ribs); (b), (c) Top-skin of wing being split into local panels and enlarged view of selected local panels; (d) Internal layout of the wing with curvilinear *SpaRibs* (concave upward); (e) Internal layout of the wing with curvilinear *SpaRibs* (convex upward); (f), (g) Top-skin of wing being split into local panels using its intersection with *SpaRibs* (concave upward) and enlarged view of selected local panels; (h), (i) Top-skin of wing being split into local panels using its intersection with *SpaRibs* (convex upward) and enlarged view of selected local panels.

### Creating Stiffened Panels and Optimization

Shell structures like the skin of a wing are usually laid with blade-stiffeners (which are similar to spars but without full depth), which help to reduce overall weight. In an aircraft wing,the stiffeners are usually placed along the span-wise direction from the wing-root to the wing-tip, as shown in Figure 11(a-c).

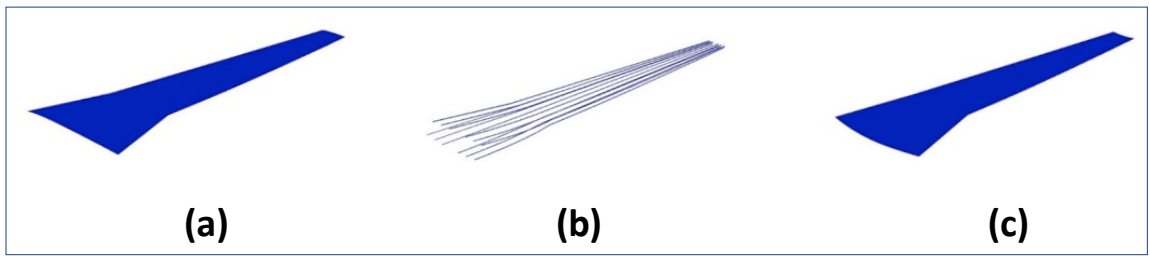

**Figure 11. (a) Upper skin (b) Blade stiffeners (c) Lower skin.**

It is challenging to determine the parts of stiffeners that are attached to each of the panels isolated for distributed optimization. An approach is developed where stiffeners on the top and bottom skin of the wing and are meshed with the quad elements. The element size is chosen to be larger than the height of the stiffeners to ensure that the stiffeners are represented by a chain of quad elements with each element sharing two common nodes with the skin. The nodes belonging to each of the local panels are already known by the method of creating local panels using the Mesh-continuity algorithm.MSC.NASTRAN*".bdf"* file for the stiffeners is generated and read to find the quad elements with nodes common with the nodes shared with the local panels.These elements form the stiffeners for respective local panels. An example of such a stiffened panel is shown in Figure 12.

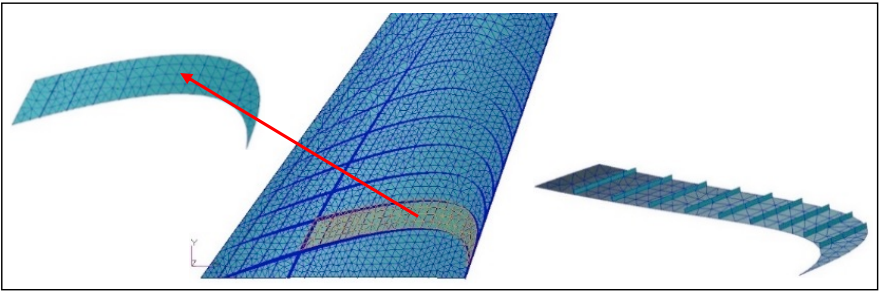

**Figure 12. Local panel with stiffeners (formed by Quad elements)**

If the cross-section of the stiffeners is rectangular, and isotropic material is used for both stiffener and the win-skin, it is reasonable to assign 3 design variables to each of the local panels: one for the thickness of the panel, one for stiffener height and one for stiffener width. In such case the total number of design, $N_{variable}$ variables is given by:

$$N_{variable} = 3(N_{Spars} + 1)(N_{Ribs} + 1) \qquad (1)$$

where, $N_{Spars}$ and $N_{Ribs}$ are the number of spars and the number of ribs, respectively.

Even though the total number of design variables can be large, in the EBF3GLWingOpt framework, three design variables associated with each panel is optimized at once. The isolated panel can be assumed to be "simply-supported," and the panel thickness, stiffener height and profile can be optimized considering constraint on stress and buckling factor. In the EBF3GLWingOpt framework, the maximum von Mises stress $(\sigma_{vm})_{max}$ and the first buckling eigenvalue $(\lambda_p)$ are computed using MSC.NASTRAN Sol 101 (linear static solver) and MSC.NASTRAN Sol 105 (buckling solver) respectively.

The optimization problem for each of the panel can be stated as follows:

***Objective:*** *min(Weight)*
***Constraint:*** $\lambda_p > 1.05$
$(\sigma_{vm})max < \sigma_y$
$(t_{stiff})_{min} \leq t \leq (t_{stiff})_{min}$
$(t_{stiff})_{min} \leq t_{stiff} \leq (t_{stiff})_{max}$
$(h_{stiff})_{min} \leq h_{stiff} \leq (h_{stiff})_{max} \qquad (2)$

where, t, $t_{stiff}$, $h_{stiff}$ are panel thickness, stiffener thickness, and stiffener height respectively and $\sigma_y$ is material yield stress.

The optimal value of t, $t_{stiff}$, $h_{stiff}$ can be obtained by gradient-based optimization or by using the approximate method as used by Jrad et al.[27]. In each iteration, the global wing model is constructed by assembling the optimized local panels and the new global deformation field is computed by running MSC.NASTRAN Sol 144 (static aeroelastic solver). From the deformation field, boundary conditions for the panels are computed for the next iteration. After every iteration, the percentage difference between the weight of the new design and the design from the previous iteration is calculated. Convergence is said to have been achieved when the difference falls within a certain user-specified limit. The optimization process is shown in Figure 13.

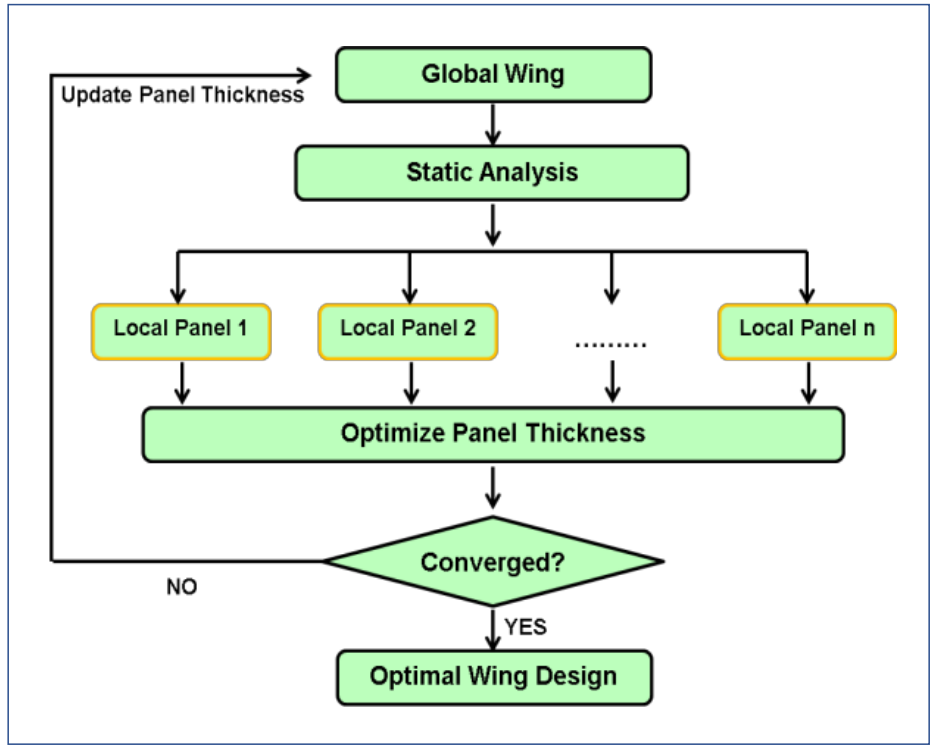

**Figure 13. Distributed design optimization of Aircraft Wing using EBF3PanelOpt (Jrad et al. [27]).**

It is worth mentioning that, the analysis of isolated stiffenened panel could further be improved by isogeometric analysis [28-33].

Although in this global-local optimization framework, the satisfaction of buckling and stress constraints is ensured for the local panel, it does not automatically lead to the satisfaction of stress and buckling constraints for the entire wing. It was found that even though weight convergence can be achieved, the global buckling constraint is often not satisfied. This drawback was overcome and the algorithm was successfully implemented by De et al. [34] to perform distributed optimization of transport aircraft wing considering 6 degrees and -2 degree angle of attack at Mach 0.85 and 35000 ft cruising altitude.

# CONCLUSIONS

This paper gives a detailed description of a new method to create stiffened local panels from two-dimensional finite element models. The capability of the algorithm used in the original EBF3GLWingOpt framework to create local panels was limited to the generation of four-edged panels with two edges along adjacent spars from the wing-skin of the CRM Wing. The new methodology is based on performing set operations on the element connectivity data of the wing's finite element model and is independent of the nodal coordinate. It can be used to create local panels of any shape and size and can be used for structures other than the NASA CRM wing. Elements of the stiffeners attached to each of the local panels are determined, and stiffened panels are created. Once the local panels are determined, the structure is ready for distributed MDO, as was shown by De et al. [35-37]. Although the algorithm has been demonstrated on the skin of a commercial aircraft wing, it is a much-generalized method and can pretty much be used to decompose a thin structure meshed with triangular shell elements. One potential application is optimization of the side frame of truck chassis with regions of variable thickness [38-41].

# ACKNOWLEDGMENTS

We are also thankful to Dr. Balakrishnan Devarajan and Varakini Sanmugadas for their suggestions to improve this article.

# DISCLOSURE STATEMENT

The authors do not report any potential conflict of interest.

# FUNDING

This work was funded by the NASA SBIR/STTR program, under contracts NNX14CD16P and NNX15CD08C.

# PUBLICATION NOTE

Part of this work was presented in the 13th World Congress on Computational Mechanics (WCCM XIII), and 2nd Pan American Congress on Computational Mechanics held in New York City from July 22 to July 27, 2018.

# REFERENCES

1. Schmit, L. A. 1960. "Structural Design by Systematic Synthesis." Pittsburg.
2. Schmit, L. A. 1981. "Structural Synthesis: Its Genesis and Development." *AIAA Journal* 19:1249–1263.
3. Schmit, L. A., and Ramanathan, R. K. 1978. "Multilevel Approach to Minimum Weight Design Including Buckling Constraints." *AIAA Journal* 16:97–104.
4. Haftka, R. T., Starnes, J. H. Jr., Barton, F. W., and Dixon, S. C. 1975. "Comparison of Two Types of Optimization Procedures for Flutter Requirements." *AIAA Journal* 13:1333–1339.
5. Starnes, J. H. Jr., and Haftka, R. T. 1979. "Preliminary Design of Composite Wings for Buckling, Strength, and Displacement Constraints." *Journal of Aircraft* 16:564–570.

6. Haftka, R. T. 1977. "Optimization of Flexible Wing Structures Subject to Strength and Induced Drag Constraints." *AIAA Journal* 15:1101–1106.
7. Aru, S., Jadhav,P., Jadhav,V., Kumar, A., and Angane P., "Design, analysis and optimization of a multi-tubular space frame." *International Journal of Mechanical and Production Engineering Research and Development (IJMPERD) ISSN (P)* (2014): 2249-6890.
8. Li, Q., Devarajan, B., Zhang, X., Burgos, R., Boroyevich, D. and Raj, P., 2016. Conceptual Design and Weight Optimization of Aircraft Power Systems with High-Peak Pulsed Power Loads (No. 2016-01-1986). SAE Technical Paper.
9. Devarajan, B., Locatelli, D., Kapania, R.K. and Meritt, R.J., 2016. Thermo-Mechanical Analysis and Design of Threaded Fasteners. In 57th AIAA/ASCE/AHS/ASC Structures, Structural Dynamics, and Materials Conference (p. 0579).
10. Singh, S., and Singh, B., "Design Optimization and Statistical Crash Analysis of Chassis Frame for off Road Vehicle" International Journal of Automobile Engineering Research and Development (IJAuERD) ISSN(P): 2277-4785
11. Fulton, R. E., Sobieszczanski, J., Storaasli, O., Landrum, E. J., and Loendor, D. 1974. "Application of Computer-Aided Aircraft Design in a Multidisciplinary Environment." *Journal of Aircraft*, 11:369–370.
12. Kroo, I., Altus, S., Braun, R., Gage, P., and Sobieski, I. 1994. "Multidisciplinary Optimization Methods for Aircraft Preliminary Design." 5th AIAA Symposium on Multidisciplinary Analysis and Optimization.
13. Manning V. M. 1999. "Large-scale design of supersonic aircraft via collaborative Optimization." PhD Thesis, Stanford University.
14. Shihan, Muhammed, et al. "Experimental Investigation and Design Optimization of Face Milling Parameters on Monel K-500 by Doe Concept." *International Journal of Mechanical and Production Engineering Research and Development (IJMPERD)* 7.4 (2017): 403-410.
15. Antoine, N. E., and Kroo, I. M. 2005. "Framework for aircraft conceptual design and environmental performance studies." *AIAA Journal* 43:2100–2109.
16. Henderson, R. P., Martins, J. and Perez, R. E. 2012. "Aircraft Conceptual Design for Optimal Environmental Performance." *Aeronautical Journal* 116:1–22.
17. Alonso, J. J., and M. R. Colonno, M. R. 2012. "Multidisciplinary optimization with applications to sonic-boom minimization." *Annual Review of Fluid Mechanics* 44:505–526.
18. Martins, J. R. R. A., and Lambe, A. B. 2013. "Multidisciplinary Design Optimization: A Survey of Architectures." *AIAA Journal* 51:2049–2075.
19. Kumar, CL., Vanaja, T., and Reddy, M., "Optimization of Parametres for Weldability Strength–an Experimental Design Approach." *International Journal of Mechanical and Production Engineering Research and Development (IJMPERD), ISSN (P)* (2018): 2249-6890.
20. Sobieszczanski-Sobieski, J., and Haftka, R. T. 1997. "Multidisciplinary Aerospace Design Optimization: Survey of Recent Developments." *Structural Optimization* 14:1–23.

21. Ciampa, P. D., and Nagel, B. 2010. "Global Local Structural Optimization of Transportation Aircraft Wings." 51st AIAA/ASME/ASCE/AHS/ASC Structures, Structural Dynamics, and Materials Conference, Orlando, Florida.
22. Yang, W., Yue, Z., Li, L., and Wang, P. 2016. "Aircraft wing structural design optimization based on automated finite element modelling and ground structure approach." *Engineering Optimization* 48:94−114.
23. Barkanov, E., Eglītis, E., Almeida, F., Bowering, M. C., and Watson, G. 2016. "Weight optimal design of lateral wing upper covers made of composite materials." *Engineering Optimization* 48:1618−1637.
24. Taminger, K. M. and Hafley, R. A. 2006. "Electron Beam Freeform Fabrication for Cost Effective Near-Net Shape Manufacturing," in VA:NATO/RTOAVT-139 specialists' meeting on cost effective manufacture via net shape processing. Amsterdam (The Netherlands): NATO Hampton.
25. Locatelli, D., Mulani, S. B., and Kapania, R. K. 2011. "Wing-Box Weight Optimization Using Curvilinear Spars and Ribs (SpaRibs)." *Journal of Aircraft* 48:1671−1684.
26. Locatelli, D., Tamijani, A. Y., Mulani, S. B., and Kapania, R. K. 2013. "Multidisciplinary Optimization of Supersonic Wing Structures Using Curvilinear Spars and Ribs (SpaRibs)." 54th AIAA/ASME/ASCE/AHS/ASC Structures, Structural Dynamics, and Materials and Co-located Conferences, Boston, Massachusetts. *AIAA* 2013−1931.
27. Jrad, M., De, S., and Kapania, R. K. 2017. "Global-Local Aeroelastic Optimization of Internal Structure of Transport Aircraft Wing." 18th AIAA/ISSMO Multidisciplinary Analysis and Optimization Conference, Denver, Colorado.
28. Devarajan, B. and Kapania, R.K., 2020. Thermal buckling of curvilinearly stiffened laminated composite plates with cutouts using isogeometric analysis. Composite Structures, 238, p.111881.
29. Miglani, J., Devarajan, B. and Kapania, R.K., 2018. Thermal buckling analysis of periodically supported composite beams using Isogeometric analysis. In 2018 AIAA/ASCE/AHS/ASC Structures, Structural Dynamics, and Materials Conference (p. 1224).
30. Miglani, J., Devarajan, B., & Kapania, R. K. (2021). Isogeometric Thermal Buckling and Sensitivity Analysis of Periodically Supported Laminated Composite Beams. *AIAA Journal*, 1-10.
31. Devarajan, B. (2021). Vibration Analysis of Timoshenko Beams using Isogeometric Analysis. *arXiv preprint arXiv:2104.12860*.
32. Devarajan, B. (2021). Free Vibration analysis of Curvilinearly Stiffened Composite plates with an arbitrarily shaped cutout using Isogeometric Analysis. *arXiv preprint arXiv:2104.12856*.
33. Devarajan, B. (2021). Analyzing Thermal Buckling in Curvilinearly Stiffened Composite Plates with Arbitrary Shaped Cutouts Using Isogeometric Level Set Method. *arXiv preprint arXiv:2104.05132*.
34. De, S., Jrad, M., and Kapania, R. K. 2019. "Structural Optimization of Internal Structure of Aircraft Wings with Curvilinear Spars and Ribs." *Journal of Aircraft* 56:707−718.
35. Robinson, J. H., Doyle, S., Ogawa, G., Baker, M., De, S., Jrad, M., & Kapania, R. K. (2016). Aeroelastic Optimization of Wing Structure Using Curvilinear Spars and Ribs (SpaRibs). In *17th AIAA/ISSMO Multidisciplinary Analysis and Optimization Conference* (p. 3994).

36. De, S., Jrad, M., Locatelli, D., Kapania, R. K., and Baker, M. 2017. "SpaRibs Geometry Parameterization for Wings with Multiple Sections using Single Design Space." 58th AIAA/ASCE/AHS/ASC Structures, Structural Dynamics, and Materials, Dallas.
37. De, S. (2017). *Structural Modeling and Optimization of Aircraft Wings having Curvilinear Spars and Ribs (SpaRibs)* (Doctoral dissertation, Virginia Tech).
38. De, S.; Singh, K.; Seo, J.; Kapania, R.K.; Ostergaard, E.; Angelini, N.; Aguero, R. Lightweight Chassis Design of Hybrid Trucks Considering Multiple Road Conditions and Constraints. *World Electr. Veh. J.* **2021**, *12*, 3.
39. De, S., Singh, K., Alanbay, B., Kapania, R. K., "Structural Optimization of Truck Front-Frame under Multiple Load Cases," in ASME International Mechanical Engineering Congress & Exposition, Pittsburgh, PA, November 9-15, 2018.
40. De, S., Singh, K., Seo, J., Kapania, R. K., "Structural Design and Optimization of Commercial Vehicles Chassis under Multiple Load Cases and Constraints," in 60th AIAA/ ASCE/AHS/ASC Structures, Structural Dynamics, and Materials Conference, San Diego, CA, 2018.
41. De, S., Singh, K., Seo, J., Kapania, R. et al., "Unconventional Truck Chassis Design with Multi-Functional Cross Members," SAE Technical Paper 2019-01-0839, 2019